\documentclass[12pt]{amsart}

\title{\large{The closed topological vertex via the Cremona transform} }
\author{Jim Bryan and Dagan Karp}
\date{\today}
\address{
Department of Mathematics\\
University of British Columbia \\
Vancouver, BC, Canada 
}
\email{jbryan@math.ubc.ca}
\email{dkarp@math.ubc.ca}

\usepackage{eepic,epic}

\usepackage{amsmath}
\usepackage{amsmath,amsthm,amsfonts,verbatim}
\newtheorem{theorem}{Theorem}

\newtheorem{lemma}[theorem]{Lemma}
\newtheorem{proposition}[theorem]{Proposition}

\newtheorem{assumption}{Assumption}

\newtheorem{rem1}[theorem]{Remark}
\newenvironment{remark}{\begin{rem1}\em}{\end{rem1}}

\newcommand{\cnums} {{\mathbb C}}          
\newcommand{\znums} {{\mathbb Z}}		

\newcommand{\im}{\operatorname{Im}}

\renewcommand{\P}{\mathbb{P}}
\newcommand{\M}{\overline{{M}}}
\renewcommand{\O}{\mathcal{O}}

\begin{document}
\pagestyle{plain}
\maketitle 

\begin{abstract}
We compute the local Gromov-Witten invariants of the ``closed vertex'',
that is, a configuration of three $\P ^{1}$'s meeting in a single triple
point in a Calabi-Yau threefold. The method is to express the local
invariants of the vertex in terms of ordinary Gromov-Witten invariants of a
certain blowup of $\P ^{3}$ and then to compute those invariants via the
geometry of the Cremona transformation.
\end{abstract}



\section{Introduction}

Let $C\subset Y$ be a curve in a Calabi-Yau threefold consisting of three
$\P ^{1}$'s meeting in a single point $p$. With some minor assumptions on
the formal neighborhood of $C\subset Y$ (see Assumption~\ref{assumption:
local structure of C}), the contribution to the genus $g$ Gromov-Witten
invariants of $Y$ by maps to $C$ is well defined and denoted by
$N_{d_{1},d_{2},d_{3}}^{g} (C)$, where $d_{i}$ is the degree of the map to
$i$th component. We call $C$ the closed topological vertex and we say that
$N_{d_{1},d_{2},d_{3}}^{g} (C)$ are the local invariants of $C$.

\begin{theorem}\label{thm: main thm}
The local invariants of the closed topological vertex are given as follows.
\[
N^{g}_{d_{1},d_{2},d_{3}} (C) = 0
\]
if $\{d_{1},d_{2},d_{3} \}$ contains two distinct non-zero values, otherwise
\[
N^{g}_{d,d,d} (C)=N^{g}_{d,d,0} (C)=N^{g}_{d,0,0} (C).
\]
\end{theorem}

Note that $N^{g}_{d,0,0} (C)$ is just the local invariant of a smooth $\P
^{1}$ in a Calabi-Yau threefold embedded with normal bundle $\O (-1)\oplus
\O (-1)$. These were computed by Faber and Pandharipande \cite{Fa-Pa} and
are given by:
\[
N_{d,0,0}^{g} (C)=\tfrac{|B_{2g} (2g-1)|}{(2g)!}d^{2g-3}
\] 
where $B_{2g}$ is the $2g$th Bernoulli number.

If one of the $d_{i}$'s is zero, then Theorem~\ref{thm: main thm} is a
special case of the computation of local invariants of ADE configurations
done by Bryan, Katz, and Leung in \cite{BKL}.

The local invariants $N^{g}_{d_{1},d_{2},d_{3}} (C)$ are related to the
``topological vertex'' of Aganagic, Klemm, Marino, and Vafa
\cite{AKMV}. Using the conjectural Chern-Simons/string theory duality, they
compute the open string amplitudes of 3 rational curves meeting in a triple
point in a Calabi-Yau threefold. Open string amplitudes correspond to a
version of Gromov-Witten theory using Riemann surfaces with boundary. The
correct mathematical formulation of open string Gromov-Witten theory is not
currently known. Our local invariants $N^{g}_{d_{1},d_{2},d_{3}} (C)$
correspond to closed string amplitudes, which is why we call our
configuration the closed topological vertex.

In the course of the proof of Theorem~\ref{thm: main thm}, we prove a
certain symmetry of the Gromov-Witten of $\P ^{3}$ blown up at points. This
symmetry arises from the Cremona transformation and may be of independent
interest.

Let $X$ be the blowup of $\P ^{3}$ at six (distinct) points $\{x_{1},\dots
,x_{6} \}$. Let $h$ be the pullback of the class of a line in $\P ^{3}$ and
let $e_{1},\dots ,e_{6}$ be the classes of the lines in the corresponding
exceptional divisors.

\begin{theorem}\label{thm: cremona invariance}
Let $ \beta =dh -\sum _{i=1}^{6}a_{i}e_{i}$ with $2d=\sum _{i=1}^{6}a_{i}$
and assume that either $a_{5}$ or $a_{6}$ is non-zero. Then we have the
following equality of Gromov-Witten invariants:
\[
\left\langle \; \right\rangle^{X}_{g,\beta }=\left\langle \;
\right\rangle^{X}_{g,\beta '}
\]
where $\beta '=d'h-\sum _{i=1}^{6}a_{i}'e_{i}$  has coefficients given by
\begin{align*}
d'\, &=3d-2 (a_{1}+a_{2}+a_{3}+a_{4})\\
a_{1}'&=\;\, d-\; \, (a_{2}+a_{3}+a_{4})\\
a_{2}'&=\;\, d- \;\, (a_{1}+a_{3}+a_{4})\\
a_{3}'&=\;\, d-\; \, (a_{1}+a_{2}+a_{4})\\
a_{4}'&=\;\, d-\; \, (a_{1}+a_{2}+a_{3})\\
a_{5}'&=\;\, a_{5}\\
a_{6}'&=\;\, a_{6}.
\end{align*}
\end{theorem}

\begin{remark}
It is not clear if the above result generalizes to more than six
points. The condition that $a_{5}$ or $a_{6}$ is non-zero is necessary; for
example, the theorem fails for the class $\beta =h-e_{1}-e_{2}$.
\end{remark}

Our basic strategy of first identifying the local invariant with an
invariant of a blowup of projective space and then utilizing the Cremona
transformation was first employed in \cite{Br-Le1}. In that paper, the
technique was used to compute the local contributions of nodal fibers in an
elliptically fibered $K3$ surface (see section~5.3 of \cite{Br-Le1})

\subsubsection*{Acknowledgments:} The authors give warm thanks to Jacob
Shapiro, Kalle Karu, and S\'andor K\'ovac for helpful conversations. The
first author thanks NSERC and the Sloan Foundation for their support, and
the second author thanks IPAM for its hospitality and support.

\section{Local invariants}\label{sec: local invariants}

Let $i:Z\subset Y$ be a closed subvariety of a smooth projective Calabi-Yau
threefold. Let $\beta $ be a curve class in $Z$; we define the ``local
Gromov-Witten invariants of $Z\subset Y$'', denoted $N_{\beta }^{g}
(Z\subset Y)$, whenever the substack of $\M _{g} (Y,i_{*}\beta )$
consisting of maps whose image lies in $Z$ is a connected component. This
component then inherits a degree 0 virtual class and $N_{\beta }^{g}
(Z\subset Y)$ is defined to be its degree. In general, $N^{g}_{\beta }
(Z\subset Y)$ depends on a (formal or analytic) neighborhood of $Z\subset
Y$, but in many cases, it depends only on the normal bundle. We write
$N^{g}_{\beta } (Z) $ when the structure of the neighborhood is understood.

Local invariants have been studied extensively in the literature. Examples
include the local invariants of surfaces such as $\P ^{2}$
\cite{Graber-Zaslow,Katz-Klemm-Vafa}, $K3$ \cite{Shapiro}, and rational
elliptic surfaces \cite{HST}. Local invariants of curves have been studied
in \cite{BKL,Br-Pa-TQFT,Br-Pa,
Pandharipande-degenerate-contributions}. Recently, Aganagic, Klemm, Marino,
and Vafa (\cite{AKMV}, c.f. \cite{Diaconescu-Florea}) have developed a
technique, based on a conjectural physical duality, which computes the
local invariants of any toric curve or surface. The basic building block in
their algorithm are the open string amplitudes of the vertex configuration
$C$. The correct mathematical definition of the open-string invariants is
not currently known.

In order to define $N^{g}_{d_{1},d_{2},d_{3}} (C)$ we must specify the
geometry of the formal neighborhood of $C$ in $Y$.

\begin{assumption}\label{assumption: local structure of C}
We assume that the local geometry of $C\subset Y$ is as follows. The curve
$C$ consists of three components $C=C_{1}\cup C_{2}\cup C_{3}$, with
$C_{i}\cong \P ^{1}$, and meeting in a single triple point $p$. $C$ is
embedded in $Y$ such that the normal bundle of each component of $C$ is
isomorphic to $\O (-1)\oplus \O (-1)$. Additionally, for the case of
$d_{i}>0$, we assume the formal neighborhood of the triple point has the
geometry of the coordinate axes in $\cnums ^{3}$ with respect to the local
coordinates defined by the normal bundles.  For the case where one of the
$d_{i}$'s is zero, say $d_{3}$, we assume that the curve $C_{1}\cup C_{2}$
is contractable (c.f. \cite{BKL}, section 2).
\end{assumption}

We remark that the above assumption leads to two different formal
neighborhoods in the two cases $d_{i}>0$ and $d_{3}=0$. In the former case,
the curve $C_{1}\cup C_{2}$ admits a deformation in $Y $ that smooths the
node and hence it cannot be contractable.

\section{The closed vertex and invariants of a blowup of $\P ^{3}$}\label{sec:
closed vertex and blowups}

In this section we prove that the local invariants of the closed
topological vertex are equal to certain ordinary Gromov-Witten invariants:
\begin{proposition}\label{prop: N(C) invariant in X}
Let $X\to \P ^{3}$ be the blowup of $\P ^{3}$ at six points
$\{x_{1},x_{2},x_{3},x_{1}',x_{2}',x_{3}' \}$. Let $h$ be the pullback of
the class of a line in $\P^{3}$ and let $\{
e_{1},e_{2},e_{3},e_{1}',e_{2}',e_{3}'\}$ be the classes of the lines in
the corresponding exceptional divisors. Assume that $d_{1},d_{2},d_{3}>0$
and let
\[
\beta =\sum _{i=1}^{3}d_{i} (h-e_{i}-e'_{i}).
\]
Then the local invariants of the vertex are given by the ordinary
Gromov-Witten invariants of $X$ in the class of $\beta $:
\[
N^{g}_{d_{1},d_{2},d_{3}} (C)=\left\langle \; \right\rangle^{X}_{g,\beta }.
\]
\end{proposition}

\textsc{Proof:} By the deformation invariance property of Gromov-Witten
invariants, we may replace $X$ with any threefold that is deformation
equivalent to the blowup of $\P ^{3}$ at six points. We will construct such
an $X$ that is compatible with the toric structure of $\P ^{3}$. Let
$x_{0},x_{1},x_{2},x_{3}$ be the fixed points of the standard $T= (\cnums
^{\times })^{3}$ action on $\P ^{3}$. Let $X'$ be the blowup of $\P ^{3}$
at the three points $\{x_{1},x_{2},x_{3} \}$, let $C''_{i}\subset \P ^{3}$
be the line between $x_{0}$ and $x_{i}$, and let $C_{i}'\subset X'$ be the
proper transform of $C_{i}''$. Let $x'_{i}$ be the intersection of $C_{i}'$
with the exceptional divisor of $X'\to \P ^{3}$. Define $X$ to be the
blowup of $X'$ at the three points $\{x'_{1},x'_{2},x'_{3} \}$ and let
$C_{i} $ be the proper transform of $C_{i}'$. Clearly, $X$ is deformation
equivalent to a blowup of $\P ^{3}$ at six distinct points.

We choose the basis $\{H,E_{1},E_{2},E_{3},E'_{1},E'_{2},E'_{3} \}$ for
$H_{4} (X;\znums )$ where $H$ is the pullback of a hyperplane in $\P ^{3}$,
$E_{i}'$ are the exceptional divisors of $X\to X'$, and $E_{i}$ are the
total transforms of the exceptional divisors of $X'\to \P ^{3}$. We choose
the basis $\{h,e_{1},e_{2},e_{3},e'_{1},e'_{2},e'_{3} \}$ for $H_{2}
(X;\znums )$ where $h$ is a line in $H$, $e_{i}'$ is a line in $E'_{i}$,
and $e_{i}$ is a line in $E_{i}$ which is disjoint from $x_{i}'$. In this
basis, we have $H\cdot H=h$, $E_{i}\cdot E_{i}=-e_{i}$, and all other
intersections between divisors are zero.  With respect to this basis, the
class of the curve $C_{i}$ is $h-e_{i}-e_{i}'$.

The $T$ action on $\P ^{3}$ lifts to a $T$ action on $X'$ and $X$ since the
centers of each blowup are invariant. Let $F$ be the union of all the $T$
invariant curves in $X$. The $T$ invariant curves and their classes are
configured as in Figure~1. The fixed points of $T$ correspond to
the vertices in the diagram.

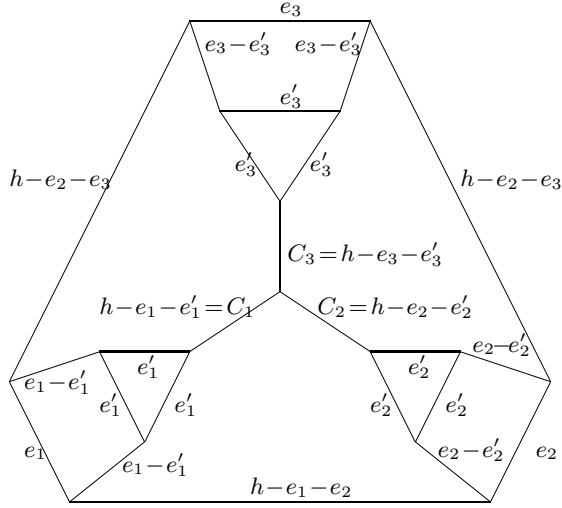
\begin{figure}\label{fig: T inv curves of X}

\setlength{\unitlength}{4mm}
\begin{picture} (20,16) (-10,0)
\put(-7 ,0 ){\line (1,0){14}}
\put(0,7 ){\line(0 ,1 ){3}}
\put(-2 ,13 ){\line(1 ,0 ){4}}
\put(-3 ,16 ){\line(1 ,0 ){6}}

\put(7 ,0 ){\line(1 ,2 ){2}}
\put(9 ,4 ){\line(-3 ,1 ){3}}
\put(6 ,5 ){\line(1 ,2 ){-1.5}}
\put(4.5,2){\line(5 ,-4 ){2.5}}
\put(3 ,5 ){\line(3,-2 ){-3}}
\put(3 ,5 ){\line(1 ,0 ){3}}
\put(0 ,10 ){\line(2 ,3 ){2}}
\put(2 ,13 ){\line(1 ,3 ){1}}
\put(3 ,16 ){\line(1 ,-2 ){6}}
\put(3 ,5 ){\line(1 ,-2 ){1.5}}

\put(-7 ,0 ){\line(-1 ,2 ){2}}
\put(-9 ,4 ){\line(3 ,1 ){3}}
\put(-6 ,5 ){\line(-1 ,2 ){-1.5}}
\put(-4.5,2){\line(5 ,4 ){-2.5}}
\put(-3 ,5 ){\line(-3,-2 ){-3}}
\put(-3 ,5 ){\line(-1 ,0 ){3}}
\put(-0 ,10 ){\line(-2 ,3 ){2}}
\put(-2 ,13 ){\line(-1 ,3 ){1}}
\put(-3 ,16 ){\line(-1 ,-2 ){6}}
\put(-3 ,5 ){\line(-1 ,-2 ){1.5}}
\scriptsize
\put(-1 ,.25 ){$h\! -\!e_{1}\! -\!e_{2}$}
\put(6,10.5 ){$h\! -\!e_{2}\! -\!e_{3}$}
\put(-9,10.5 ){$h\! -\!e_{2}\! -\!e_{3}$}
\put(0.25 ,8 ){$C_{3}\!=\!h\! -\!e_{3}\! -\!e_{3}'$}
\put(1.25 ,6.25 ){$C_{2}\!=\!h\! -\!e_{2}\! -\!e_{2}'$}
\put(-6 ,6.25 ){$h\! -\!e_{1}\! -\!e_{1}'\!=\! C_{1}$}
\put(-3.5 ,3 ){$e_{1}'$}
\put(-6 ,3 ){$e_{1}'$}
\put(-4.75 ,4.25 ){$e_{1}'$}
\put(-8.5 ,1.5 ){$e_{1}$}
\put(-5.25 ,1 ){$e_{1}\! -\!e_{1}'$}
\put(-8.5 ,3.75 ){$e_{1}\! -\!e_{1}'$}
\put(3 ,3 ){$e_{2}'$}
\put(5.5 ,3 ){$e_{2}'$}
\put(4.25 ,4.25 ){$e_{2}'$}
\put(8.5 ,1.5 ){$e_{2}$}
\put(5.25 ,1.5 ){$e_{2}\! -\!e_{2}'$}
\put(6.4 ,5 ){$e_{2}\!\! -\!\! e_{2}'$}
\put(1 ,11 ){$e_{3}'$}
\put(-1.5 ,11 ){$e_{3}'$}
\put(0 ,13.25 ){$e_{3}'$}
\put(.5 ,15 ){$e_{3}\! -\!e_{3}'$}
\put(-2.5 ,15 ){$e_{3}\! -\!e_{3}'$}
\put(0 ,16.25 ){$e_{3}$}

\end{picture}

\caption{The curve $F$, which is the union of the $T$ invariant curves of
$X$}
\end{figure}

The normal bundle of each $C_{i}$ is isomorphic to $\O (-1)\oplus \O
(-1)$. This can be seen as follows. Let $P''\subset \P ^{3}$ be a plane
containing the line $C_{i}''$ and let $P'\subset X'$ and $P\subset X$
denote the successive proper transforms.  Then $C_{i}\subset P$ so
$N_{C_{i}/P}$ defines a sub line bundle in $N_{C_{i}/X}$ of degree
$C_{i}\cdot C_{i}$ where the intersection product is in $P$. By
functoriality of blowups, $P'\to P''$ is the blowup along $x_{i}$ and
$P\to P'$ is the blowup along $x_{i}'$. We then compute the intersection
product in $P$:
\[
C_{i}\cdot C_{i}=(h-e_{i}-e_{i}')\cdot (h-e_{i}-e_{i}') =-1.
\]
Since we can apply this argument to any plane $P$ containing $C_{i}$, and
they span the normal bundle, we conclude that $N_{C_{i}/X}\cong \O
(-1)\oplus \O (-1)$. Moreover, this argument shows that the configuration
$C=C_{1}\cup C_{2}\cup C_{3}$ satisfies (the $d_{i}>0$ case of)
Assumption~\ref{assumption: local structure of C} since we can take the two
spanning planes of the normal bundle of (say) $C_{1}$ to contain $C_{2}$
and $C_{3}$ respectively and so the local geometry of the triple point is
that of the coordinate axes.

We've shown that $C\subset X$ has the local geometry of the closed vertex;
then to prove the proposition, we need to show that the only contributions
to $\left\langle \; \right\rangle^{X}_{\beta }$ are given by maps to the
configuration $C$. This is accomplished with the following:

\begin{lemma}\label{lem: all curves in beta have image in C}
Let $X$ and $C=C_{1}\cup C_{2}\cup C_{3}$ be as above and let 
\[
\beta =\sum
_{i=1}^{3} d_{i}(h-e_{i}-e_{i}')
\]
where we assume that $d_{i}>0$. Then every stable map $[f]\in \M _{g}
(X,\beta )$ has image $C$.
\end{lemma}

\textsc{Proof:}
Since the torus $T$ acts on $X$ it acts on $\M _{g} (X,\beta )$. Let 
\[
[f:D\to X]\in \M _{g} (X,\beta )
\]
be any stable map. Assume that $\im (f)$ contains a point $x$
not contained in $C$.

We study the limits of $[f]$ under the $T$ action on $\M _{g} (X,\beta
)$. For any one parameter subgroup $\cnums ^{\times }\to T$, the limit of
$x$, under the action of $\lambda \in \cnums ^{\times }$ as $\lambda \to
0$, is in the fixed point set of the $\cnums ^{\times }$ action. Moreover,
since $x\not\in C$, we can find a subgroup $\cnums ^{\times }\to T$ such
that the limit of $x$ is $v_{0}$, a fixed point of the $T$ action but
\emph{not} contained in $C$. It follows that the limit of the same $\cnums
^{\times }$ subgroup acting on $[f]\in \M _{g} (X,\beta )$ is a stable map
$f'$ such that $v_{0}\in \im (f')$. Consequently, $v_{0}$ is in the image
of all maps in the closure of the $T$ orbit of $[f']$. In particular, there
exists a stable map $[f'']\in \M _{g} (X,\beta )$ fixed by $T$ and such
that $v_{0}\in \im (f'')$. Therefore, we have constructed a stable map
$f'':D\to X$ whose image is contained in $F$, the union of the $T$
invariants curves, but is not contained in $C$.

We now show that this leads to a contradiction. The class
\[
f''_{*}[D]=\beta =\sum _{i=1}^{3}d_{i} (h-e_{i}-e_{i'})
\]
has the property that the total multiplicity of the $e$ terms is $-2
(d_{1}+d_{2}+d_{3})$ while the multiplicity of $h$ is
$(d_{1}+d_{2}+d_{3})$. Observe that every component of $F$ whose class has
an $h$ also has two negative $e$ terms. It follows that the sum of all the
non-$h$ components in $f''_{*}[D]$ must have the same number of positive
$e$ terms as negative $e$ terms. Since the classes of the non-$h$
components are all either $e_{i}$, $e_{i}'$, or $e_{i}-e_{i}'$, we can
conclude that \emph{ $\im (f'')\subset F$ does not contain any of the
components of class $e_{i}$ or $e_{i}'$.}

\begin{figure}\label{fig: T inv curves of X minus e and e' terms}
\setlength{\unitlength}{4mm}
\begin{picture} (20,16) (-10,0)
\put(-7 ,0 ){\line (1,0){14}}
\put(0,7 ){\line(0 ,1 ){3}}
\put(-3 ,5 ){\line(-3,-2 ){-3}}
\put(9 ,4 ){\line(-3 ,1 ){3}}
\put(4.5,2){\line(5 ,-4 ){2.5}}
\put(3 ,5 ){\line(3,-2 ){-3}}
\put(2 ,13 ){\line(1 ,3 ){1}}
\put(3 ,16 ){\line(1 ,-2 ){6}}
\put(-9 ,4 ){\line(3 ,1 ){3}}
\put(-4.5,2){\line(5 ,4 ){-2.5}}
\put(-2 ,13 ){\line(-1 ,3 ){1}}
\put(-3 ,16 ){\line(-1 ,-2 ){6}}
\scriptsize
\put(-1 ,.25 ){$h\! -\!e_{1}\! -\!e_{2}$}
\put(6,10.5 ){$h\! -\!e_{2}\! -\!e_{3}$}
\put(-9,10.5 ){$h\! -\!e_{2}\! -\!e_{3}$}
\put(0.25 ,8 ){$C_{3}$}
\put(1.25 ,6.25 ){$C_{2}$}
\put(-2.25 ,6.25 ){$C_{1}$}
\put(-5.25 ,1 ){$e_{1}\! -\!e_{1}'$}
\put(-8.5 ,3.75 ){$e_{1}\! -\!e_{1}'$}
\put(5.25 ,1.5 ){$e_{2}\! -\!e_{2}'$}
\put(6.4 ,5 ){$e_{2}\!\! -\!\! e_{2}'$}
\put(.5 ,15 ){$e_{3}\! -\!e_{3}'$}
\put(-2.5 ,15 ){$e_{3}\! -\!e_{3}'$}
\end{picture}
\caption{Possible curves in $\im (f'')$}
\end{figure}
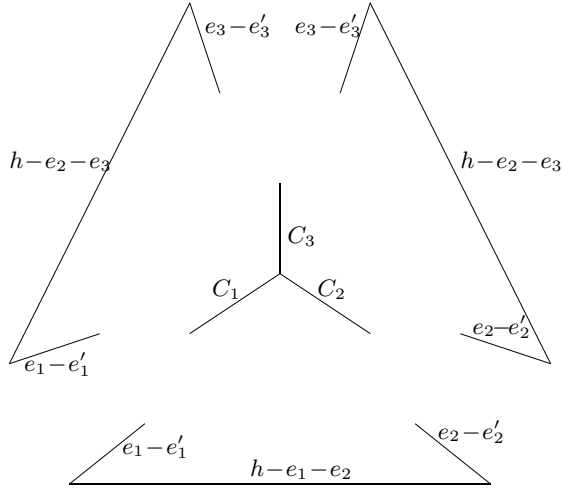

Removing those components from the diagram of $F$, we find that $\im (f'')$
must lie in the components pictured in Figure~2. Since the domain of $f''$
is connected, the image is connected and must lie in one of the components
of the graph shown in Figure~2. Since $\im (f'')\not\subset C$, it must be
one of the outer components. However, this contradicts the assumption that
$d_{i}>0$ for all $i$ since none of the outer components contain all three
classes $\{e_{1},e_{2},e_{3} \}$.

This proves that our initial assumption is false and so $\im (f)\subset C$
for all $[f]\in \M _{g} (X,\beta )$ and the lemma is proved. \qed

\section{Computing the local invariants $N^{g}_{d_{1},d_{2},d_{3}}
(C)$}\label{sec: computing N}

We need the following vanishing lemma (c.f. \cite{Gathmann}
Proposition~3.1):

\begin{lemma}\label{lem: vanishing}
Let $X$ and $\beta $ be as in Theorem~\ref{thm: cremona invariance}, then
if $a_{i}<0$ for some $i$, $\left\langle \; \right\rangle^{X}_{g,\beta
}=0$.
\end{lemma}

\textsc{Proof:} Without loss of generality we may assume that $i=1$ so that
$a_{1}<0$. Let $[f]\in \M _{g} (X,\beta )$ be any stable map. Since $\beta
\cdot E_{1}<0$, $\im (f)$ must have a component $C''$ contained in
$E_{1}$. The class of $C''$ is therefore $ne_{1}$ for some $n>0$. Writing
$\im (f)=C'\cup C''$ we have
\[
0=-K_{X}\cdot \beta =-K_{X}\cdot C'-K_{X}\cdot C''=-K_{X}\cdot C'+2n
\]
and so $-K_{X}\cdot C'$ is negative. However, this contradicts
Lemma~\ref{lem: -K is nef on X} which states that $-K_{X}$ is
nef. Therefore $\M _{g} (X,\beta )=\emptyset $ and the Lemma follows.\qed

We now prove Theorem~\ref{thm: main thm} assuming Theorem~\ref{thm: cremona
invariance} which we will prove in the subsequent sections.

We can assume that $d_{1},d_{2},d_{3}>0$ since the other cases are covered
by \cite{BKL} and \cite{Fa-Pa}. We order the $d_{i}$'s so that $d_{1}\geq
d_{2}\geq d_{3}$.

By Proposition~\ref{prop: N(C) invariant in X} we have
\[
N^{g}_{d_{1},d_{2},d_{3}} (C)=\left\langle \;  \right\rangle^{X}_{g,\beta }
\]
where
\[
\beta = (d_{1}+d_{2}+d_{3})h-d_{1} (e_{1}+e_{2})-d_{2} (e_{3}+e_{4})-d_{3}
(e_{5}+e_{6}).
\]
Then applying Theorem~\ref{thm: cremona invariance} we get 
\[
N^{g}_{d_{1},d_{2},d_{3}} (C)=\left\langle \;  \right\rangle^{X}_{g,\beta '}
\]
where
\[
\beta '= (3d_{3}-d_{1}-d_{2})h- (d_{3}-d_{2})e_{1}- (d_{3}-d_{2})e_{2}-
(d_{3}-d_{1})e_{3}- (d_{3}-d_{1})e_{4}-d_{3}e_{5}-d_{3}e_{6}.
\]
Lemma~\ref{lem: vanishing} then implies that $N^{g}_{d_{1},d_{2},d_{3}}
(C)=0$ unless $d_{3}=d_{2}=d_{1}=d$, in which case
\[
\beta '=d (h-e_{5}-e_{6})
\]
and so $N^{g}_{d,d,d} (C)=N^{g}_{0,0,d} (C)$ and Theorem~\ref{thm: main
thm} is proved. \qed

\section{The geometry of the Cremona transformation}\label{sec: geometry of
Cremona}

In this section, we prove Theorem~\ref{thm: cremona invariance} by studying
the geometry of $X$, the blowup of $\P ^{3}$ at six points, and $\hat{X}$,
the blowup of $X$ along a certain configuration of six lines.

Let $X$ be the blowup of $\P ^{3}$ at six distinct points $x_{1},\dots
,x_{6}$. We take the first four points to be the fixed points of the
standard torus action on $\P ^{3}$ and we take the remaining two points to
be any fixed points of the Cremona transformation:

\begin{align*}
\P ^{3}& \dasharrow  \P ^{3}\\
(z_{0}:z_{1}:z_{2}:z_{3})&\mapsto
(\frac{1}{z_{0}}:\frac{1}{z_{1}}:\frac{1}{z_{2}}:\frac{1}{z_{3}}).
\end{align*}
Let $l_{jk}$, $1\leq j<k\leq 4$ be the proper transform of the line through
$x_{j}$ and $x_{k}$. Let
\[
\pi :\hat{X}\to X 
\]
be the blowup of $X$ along the six (disjoint)
lines $l_{jk}$.

$\hat{X}$ admits an involution $\tau :\hat{X}\to \hat{X}$ which resolves
the Cremona transformation.  The map $\tau $ is discussed in more detail by
Gathmann in \cite{Gathmann}, although note that our $\hat{X}$ has the two
additional blowups at $x_{5}$ and $x_{6}$ whose corresponding exceptional
divisors are simply fixed by $\tau $.

We briefly describe the divisors and the curves on $X$ and $\hat{X}$ and
their intersections. Generally we denote divisor classes with upper case
letters and curve classes with lower case letters. Classes on $\hat{X}$
will have a hat, and classes on $X$ will not.

The homology groups $H_{4} (X;\znums )$ and $H_{2} (X;\znums )$ are spanned
by the divisor and curve classes respectively:
\[
H_{4} (X;\znums )=\left\langle H,E_{1},\dots ,E_{6} \right\rangle, \quad 
H_{2} (X;\znums )=\left\langle h,e_{1},\dots ,e_{6} \right\rangle.
\]
Here $H$ is the pullback of the hyperplane in $\P ^{3}$, $h$ is the
class of the line in $H$, $E_{i}$ is the exceptional divisor over $x_{i}$,
and $e_{i}$ is the class of a line in $E_{i}$.

The intersection pairing on $X$ is given by:
\begin{align*}
H\cdot H&=h,& E_{i}\cdot E_{i}&=-e_{i},\\
H\cdot h&=p,& E_{i}\cdot e_{i}&=-p
\end{align*}
where $p\in H_{0} (X;\znums )$ is the class of the point and all other
pairings are zero.

The homology groups $H_{4} (\hat{X};\znums )$ and $H_{2} (\hat{X};\znums )$
are also spanned by divisor and curve classes:
\[
H_{4} (\hat{X};\znums )=\left\langle \hat{H},\hat{E}_{i},\hat{F}_{jk}
\right\rangle,\quad H_{2}
(\hat{X};\znums )=\left\langle \hat{h},\hat{e}_{i},\hat{f}_{jk}
\right\rangle,
\]
where $1\leq i\leq 6$ and $1\leq j<k\leq 4$.  Here $\hat{H}$ is the proper
transform of $H$ and $\hat{h}$ is the generic line in
$\hat{H}$. $\hat{E}_{i}$ is the proper transform of $E_{i}$ and
$\hat{e}_{i}$ is the class of the generic line in $\hat{E}_{i}$.
$\hat{F}_{jk}$ is the component of the exceptional divisor of $\hat{X}\to
X$ lying over $l_{jk}$, and $\hat{f}_{jk}$ is the fiber class of $\pi
:\hat{F}_{jk}\to l_{jk}$.

Note that $\hat{F}_{jk}\to l_{jk}$ is the trivial fibration and the class
of the section $\hat{s}_{jk}$ is given by
\[
\hat{s}_{jk}=\hat{h}-\hat{e}_{j}-\hat{e}_{k}+\hat{f}_{jk}.
\]

The intersections are given as follows:
\begin{align*}
\hat{H}\cdot \hat{H}\; \;  &=\hat{h},& \hat{E}_{i}\cdot
\hat{E}_{i}\; \, &=-\hat{e_{i}},&
\hat{F}_{jk}\cdot \hat{F}_{jk}&=-\hat{s}_{jk}-\hat{f}_{jk},\\
\hat{H}\cdot \hat{F}_{jk}&=\hat{f}_{jk},& \hat{E}_{j}\cdot
\hat{F}_{jk}&=\; \;  \hat{f}_{jk},&&\\
\hat{H}\cdot \hat{h}\; \; \; \, &=\hat{p},&\hat{E}_{i}\cdot \hat{e}_{i}\;
\; \, &=-\hat{p},&\hat{F}_{jk}\cdot \hat{f}_{jk}\; &=-\hat{p},
\end{align*}
where $\hat{p}\in H_{0} (\hat{X};\znums )$ is the class of the point and
all other intersections are zero.

The action of $\tau $ on the curve classes of $\hat{X}$ is given by:
\begin{align*}
\tau _{*}\hat{h}\; &=3\hat{h}+2 (\hat{e}_{1}+ \hat{e}_{2}+\hat{e}_{3} +\hat{e}_{4}),\\
\tau _{*}\hat{e}_{1}&=\;\, \hat{h}+\; \, (\hat{e}_{2}+\hat{e}_{3} +\hat{e}_{4}),\\
\tau _{*}\hat{e}_{2}&=\;\, \hat{h}+\;\,  (\hat{e}_{1}+\hat{e}_{3} +\hat{e}_{4}),\\
\tau _{*}\hat{e}_{3}&=\;\, \hat{h}+\;\,  (\hat{e}_{1}+\hat{e}_{2} +\hat{e}_{4}),\\
\tau _{*}\hat{e}_{4}&=\;\, \hat{h}+\;\,  (\hat{e}_{1}+\hat{e}_{2} +\hat{e}_{3}),\\
\tau _{*}\hat{e}_{5}&=\; \,\hat{e}_{5},\\
\tau _{*}\hat{e}_{6}&=\; \,\hat{e}_{6},\\
\tau _{*}\hat{f}_{jk}&=\; \,\hat{s}_{jk}.
\end{align*}

For a class $\hat{\beta }=d\hat{h}-\sum _{i=1}^{6}a_{i}\hat{e}_{i}$ with
$2d=\sum _{i=1}^{6}a_{i}$, we have $-K_{\hat{X}}\cdot \hat{\beta }=0$ and
so the degree $\hat{\beta }$ Gromov-Witten invariants have no insertions.
Since $\tau $ is an isomorphism, it preserves the Gromov-Witten invariants
of $\hat{X}$ so in particular, 

\[
\left\langle\;  \right\rangle^{\hat{X}}_{g,\hat{\beta }}=\left\langle
\; \right\rangle^{\hat{X}}_{g,\tau _{*}\hat{\beta }}
\]
where 
\[
\tau _{*}\hat{\beta }=d'\hat{h}-\sum _{i=1}^{6}a_{i}'\hat{e}_{i}
\]
has
coefficients $d',a_{1}',\dots ,a_{6}'\; $ given by the equations of
Theorem~\ref{thm: cremona invariance}.

To prove Theorem~\ref{thm: cremona invariance} then, it suffices to prove
the following
\begin{lemma}\label{lem: invariant on Xhat = inv on X}
Let $d,a_{1},\dots ,a_{6}$ be such that $2d=\sum _{i=1}^{6}a_{i}$ and
either $a_{5}$ or $a_{6}$ is non-zero. Then
\[
\left\langle\;  \right\rangle^{X}_{g,\beta }=\left\langle
\; \right\rangle^{\hat{X}}_{g,\hat{\beta }}
\]
where $\beta =dh-\sum _{i=1}^{6}a_{i}e_{i}$ and $\hat{\beta }=d\hat{h}-\sum
_{i=1}^{6}a_{i}\hat{e}_{i}$.
\end{lemma}

\begin{remark}
The condition that $a_{5}$ or $a_{6}$ is non-zero is necessary. For
example\footnote{This example provides a counterexample to the main theorem
of \cite{Hu}.},
\[
1=\left\langle \;
\right\rangle^{X}_{0,h-e_{1}-e_{2}}\neq \left\langle \;
\right\rangle^{\hat{X}}_{0,\hat{h}-\hat{e}_{1}-\hat{e}_{2}}=0.
\]
\end{remark}

\textsc{Proof:} Without loss of generality we may assume that $a_{5}\neq
0$. We will show that any $[\hat{f}]\in \M _{g} (\hat{X},\hat{\beta })$ has
an image which is disjoint from $\hat{F}=\cup _{j<k}\hat{F}_{jk}$, and any
$[f]\in \M _{g} (X,\beta )$ has an image which is disjoint from $l=\cup
_{j<k}l_{jk}$. It follows that the natural map $\M _{g} (\hat{X},\hat{\beta
})\to \M _{g} (X,\beta )$ induced by $\pi $ is an isomorphism of the moduli
spaces and their virtual fundamental classes. Indeed, if both $\im
(\hat{f})\cap \hat{F}=\emptyset $ and $\im (f)\cap l=\emptyset $ for all
stable maps $[\hat{f}]\in \M _{g} (\hat{X},\hat{\beta })$ and $[f]\in \M
_{g} (X,\beta )$, then both $\M _{g} (\hat{X},\hat{\beta })$ and $\M _{g}
(X,\beta )$ are canonically identified with $\M _{g} (\hat{X}\backslash
\hat{F},\hat{\beta })$.

Let $[f:C\to X]\in \M _{g} (X,\beta )$ and suppose that $\im (f)\cap
l_{jk}\neq \emptyset $ for some $j$ and $k$. $\im (f)\not\subset l_{jk}$
since $a_{5}\neq 0$ and so 
\[
f_{*} (C)=C'+nl_{jk}
\]
where $C'$ meets $l_{jk}$ in a finite set of points ($n$ can be zero
here). Let $\hat{C}'$ be the proper transform of $C'$. Since $C'\cap
l_{jk}\neq \emptyset $, we have $\hat{C}'\cdot \hat{F}_{jk}=m>0$. Therefore
we have
\[
\hat{C}'=d\hat{h}-\sum _{i=1}^{6}a_{i}\hat{e}_{i}-n
(\hat{h}-\hat{e}_{j}-\hat{e_{k}})-m\hat{f}_{jk}.
\]
Define $\{j',k' \}$ by the condition $\{j',k' \}\cup \{j,k \}=\{1,2,3,4 \}$
and let
\[
\hat{D}_{jk}=2\hat{H}- (\hat{E}_{1}+\dots
+\hat{E}_{6})-\hat{F}_{jk}-\hat{F}_{j'k'}.
\]
Then 
\[
\hat{D}_{jk}\cdot \hat{C}'=-m<0.
\]
However, this contradicts Lemma~\ref{lem: Djk is nef} which states that
$\hat{D}_{jk}$ is nef. Thus $\im (f)\cap l=\emptyset $ for all $[f]\in \M
_{g} (X,\beta )$.

We argue in a similar fashion for $\M _{g} (\hat{X},\hat{\beta })$.  Let
$[\hat{f}:C\to \hat{X}]\in \M _{g} (\hat{X},\hat{\beta })$ and suppose that
$\im (\hat{f})\cap \hat{F}_{jk}\neq \emptyset $ for some $j$ and $k$. Since
$\hat{\beta }\cdot \hat{F}_{jk}=0$, $\hat{f}_{*} (C)$ must have a component
$C''$ contained in $\hat{F}_{jk}$. We then have
\[
\hat{\beta}=\hat{f}_{*} (C) =C'+ C''
\]
where $C'$ is non-empty since $\hat{\beta }\cdot \hat{E}_{5}=a_{5}>0$.

Since $C''\subset \hat{F}_{jk}$ is an effective class in $\hat{F}_{jk}\cong
\P ^{1}\times \P ^{1}$, it is of the form
$n\hat{s}_{jk}+m\hat{f}_{jk}$ with $n,m\geq 0$ and $n+m>0$.

Define $\hat{D}_{jk}$ as above.  Then $ \hat{D}_{jk}\cdot \hat{\beta }=0$
and $\hat{D}_{jk}\cdot C''=n+m>0$ and so $\hat{D}_{jk}\cdot C'<0$,
contradicting the fact that $\hat{D}_{jk}$ is nef.

This proves that $\im (\hat{f})\cap \hat{F}=\emptyset $ for all
$[\hat{f}]\in \M _{g} (\hat{X},\hat{\beta })$ and 
Lemma~\ref{lem: invariant on Xhat = inv on X} is proved. 

This then completes the proof of Theorem~\ref{thm: cremona invariance}.

\section{Nef divisors on $X$ and $\hat{X}$}\label{sec: nef divisors} In
this section we prove the results about nef divisors on $X$ and $\hat{X}$
that were used earlier.

\begin{lemma}\label{lem: -K is nef on X}
$-K_{X}$ is a nef divisor on $X$.
\end{lemma}
\textsc{Proof\footnote{The idea for the proof of this lemma was suggested
to us by S\'andor Kov\'acs.}:} The anti-canonical divisor on $X$ is given
by
\[
-K_{X}=4H- 2(E_{1}+\dots +E_{6})
\]
and so $-K_{X}$ is nef if and only if 
\begin{align*}
D&=D'+D''\\
&= (H-E_{1}-E_{2}-E_{3})+ (H-E_{4}-E_{5}-E_{6})
\end{align*}
is nef. $D'$ and $D''$ are the proper transforms of the planes through
$\{x_{1},x_{2},x_{3} \}$ and $\{x_{4},x_{5},x_{6} \}$ respectively. Thus to
see that $D$ is nef, it suffices to check that $D\cdot C\geq 0$ for a curve
$C\subset D'$.

The class of $C$ is given by 
\[
dh-a_{1}e_{1}-a_{2}e_{2}-a_{3}e_{3}
\]
for some $d$, $a_{1}$, $a_{2}$, and $a_{3}$.

The first Chern class of the normal bundle of $D'\subset X$ is 
\[
D'\cdot D'= h-e_{1}-e_{2}-e_{3}
\]
so we get that 
\begin{align*}
D'\cdot C& = (h-e_{1}-e_{2}-e_{3})\cdot (dh-a_{1}e_{1}-a_{2}e_{2}-a_{3}e_{3})\\
&=d-a_{1}-a_{2}-a_{3}
\end{align*}
where the intersection product on the right hand side is in $D'$ which is
isomorphic to the blowup of $\P ^{2} $ at three points. The class
\[
2h-e_{1}-e_{2}-e_{3}
\]
is nef in the surface $D'$ since it is represented by an irreducible,
effective curve of non-negative square. Intersecting this class with $C$ in
$D'$ we conclude that $2d\geq a_{1}+a_{2}+a_{3}$.

Computing intersections on $X$ we get:
\begin{align*}
D\cdot C&=D'\cdot C+D''\cdot C\\
&= (d-a_{1}-a_{2}-a_{3})+d\\
&\geq 0
\end{align*}
which proves the lemma. 

\begin{remark}
A similar argument can be used to show that the anticanonical divisor is
nef for the blowup of $\P ^{3} $ at up to eight points. This is the optimal
result since $(-K_{X})^{3}$ is negative for blowups of more than eight
points.
\end{remark}

\begin{lemma}\label{lem: Djk is nef}
Let $1\leq j<k\leq 4$ and define $j',k'$ by the condition $\{j,k \}\cup
\{j',k' \}=\{1,2,3,4 \}$. Then the divisor
\[
\hat{D}_{jk}=2\hat{H}- (\hat{E}_{1}+\dots
+\hat{E}_{6})-\hat{F}_{jk}-\hat{F}_{j'k'}
\]
is nef in $\hat{X}$.
\end{lemma}

\textsc{Proof:} Let $\hat{D}'$ and $\hat{D}''$ be the proper transforms of
the planes through $\{x_{j},x_{k},x_{5} \}$ and $\{x_{j'},x_{k'},x_{6} \}$
respectively. Then
\begin{align*}
\hat{D}'&= \hat{H}- \hat{E}_{j}- \hat{E}_{k}- \hat{E}_{5}- \hat{F}_{jk}\\
\hat{D}''&= \hat{H}- \hat{E}_{j'}- \hat{E}_{k'}- \hat{E}_{6}- \hat{F}_{j'k'}
\end{align*}
so $\hat{D}_{jk}=\hat{D}'+\hat{D}''$.

To see that $\hat{D}_{jk} $ is nef, it suffices to check that
$\hat{D}_{jk}\cdot C\geq 0$ for any curve $C\subset \hat{D}'$. 

$\hat{D}'$ is isomorphic to the blowup of $\P ^{2}$ at three points. Under
this identification, the classes of the line and the three exceptional
divisors are

\[
h'=\hat{h}-\hat{f}_{jk},\quad e_{j}'=\hat{e}_{j}-\hat{f}_{jk},\quad
e_{k}'=\hat{e}_{k}-\hat{f}_{jk},\quad e_{5}'=\hat{e}_{5}.
\]

The curve $C\subset \hat{D}'$ has class 
\[
dh'-a_{j}e_{j}'-a_{k}e_{k}'-a_{5}e_{5}'
\]
and since $h'-e_{5}'$ is a nef divisor in $\hat{D}'$, we have 
\[
d\geq a_{5}.
\]

The first Chern class of the normal bundle of $\hat{D}'\subset \hat{X}$ is 
\[
(\hat{H}- \hat{E}_{j}- \hat{E}_{k}- \hat{E}_{5}-
\hat{F}_{jk})^{2}=-\hat{e}_{5} = -e_{5}'
\]
and so 
\[
\hat{D}'\cdot C=-e'_{5}\cdot (dh'-a_{j}e_{j}'-a_{k}e_{k}'-a_{5}e_{5}')=-a_{5}
\]
where the intersection product on the right hand side is on $\hat{D}'$.

Therefore
\begin{align*}
\hat{D}_{jk}\cdot C&=\hat{D}'\cdot C+\hat{D}''\cdot C\\
&=-a_{5}+d\\
&\geq 0.
\end{align*}
\qed


\end{document}